\documentclass[11 pt]{article}
\usepackage{amsmath}
\usepackage{amsthm}
\usepackage{bbm,amssymb}
\usepackage{bbold}
\usepackage{amsfonts,ifthen}
\usepackage[pdftex]{graphicx,color}
\usepackage[hyperfootnotes=false]{hyperref}

\newcommand{\floor}[1]{\lfloor {#1} \rfloor}
\newcommand{\ceil}[1]{\lceil {#1} \rceil}

\newtheorem{theorem}{Theorem}[section]
\newtheorem{prop}[theorem]{Proposition}
\newtheorem{lemma}[theorem]{Lemma}

\theoremstyle{remark}
\newtheorem*{remark}{Remark}

\theoremstyle{definition}

\def\div{{\text{div}\hspace{0.5ex}}}

\def\Diam{\,\mbox{Diam}}

\title{Strong Spherical Asymptotics for Rotor-Router Aggregation and the Divisible Sandpile}
\author{Lionel Levine\footnote{supported by an NSF Graduate Research Fellowship, and NSF 
grant DMS-0605166; {\tt levine(at)math.berkeley.edu}}~~and Yuval Peres\footnote{partially supported by NSF grant DMS-0605166; {\tt peres(at)stat.berkeley.edu}} \\ University of California, Berkeley and Microsoft Research}
\date{October 7, 2008}

\DeclareSymbolFont{AMSb}{U}{msb}{m}{n}
\DeclareMathSymbol{\C}{\mathbin}{AMSb}{"43} 
\DeclareMathSymbol{\EE}{\mathbin}{AMSb}{"45} 
\DeclareMathSymbol{\N}{\mathbin}{AMSb}{"4E} 
\DeclareMathSymbol{\PP}{\mathbin}{AMSb}{"50} 
\DeclareMathSymbol{\Q}{\mathbin}{AMSb}{"51} 
\DeclareMathSymbol{\R}{\mathbin}{AMSb}{"52} 
\DeclareMathSymbol{\Z}{\mathbin}{AMSb}{"5A}

\begin{document}

\maketitle
\renewcommand{\thefootnote}{}
\renewcommand{\thefootnote}{\arabic{footnote}}

\begin{abstract}

The rotor-router model is a deterministic analogue of random walk.  It can be used to define a deterministic growth model analogous to internal DLA.  We prove that the asymptotic shape of this model is a Euclidean ball, in a sense which is stronger than our earlier work \cite{LP}.  For the shape consisting of $n=\omega_d r^d$ sites, where $\omega_d$ is the volume of the unit ball in $\R^d$, we show that the inradius of the set of occupied sites is at least $r-O(\log r)$, while the outradius is at most $r+O(r^\alpha)$ for any $\alpha > 1-1/d$.
For a related model, the {\it divisible sandpile},
we show that the domain of occupied sites is a Euclidean ball with error in the radius a constant independent of the total mass.  For the classical abelian sandpile model in two dimensions, with $n=\pi r^2$ particles, we show that the inradius is at least $r/\sqrt{3}$, and the outradius is at most $(r+o(r))/\sqrt{2}$.  This improves on bounds of Le Borgne and Rossin.  Similar bounds apply in higher dimensions, improving on bounds of Fey and Redig.
\end{abstract}

\section{Introduction}

Rotor-router walk is a deterministic analogue of random walk, first studied by Priezzhev et al.\ \cite{PDDK} under the name ``Eulerian walkers.''  At each site in the integer lattice $\Z^2$ is a {\it rotor} pointing north, south, east or west.  A particle starts at the origin; during each time step, the rotor at the particle's current location is rotated clockwise by $90$ degrees, and the particle takes a step in the direction of the newly rotated rotor.  In rotor-router aggregation, introduced by Jim Propp, we start with $n$ particles at the origin; each particle in turn performs rotor-router walk until it reaches a site not occupied by any other particles.   Let $A_n$ denote the resulting region of $n$ occupied sites.  For example, if all rotors initially point north, the sequence will begin $A_1 = \{(0,0)\}$, $A_2 = \{(0,0),(1,0)\}$, $A_3 = \{(0,0),(1,0),(0,-1)\}$.  The region $A_{1,000,000}$ is pictured in Figure~1.  In higher dimensions, the model can be defined analogously by repeatedly cycling the rotors through an ordering of the $2d$ cardinal directions in $\Z^d$.

\begin{figure}
\centering
\includegraphics[scale=.25]{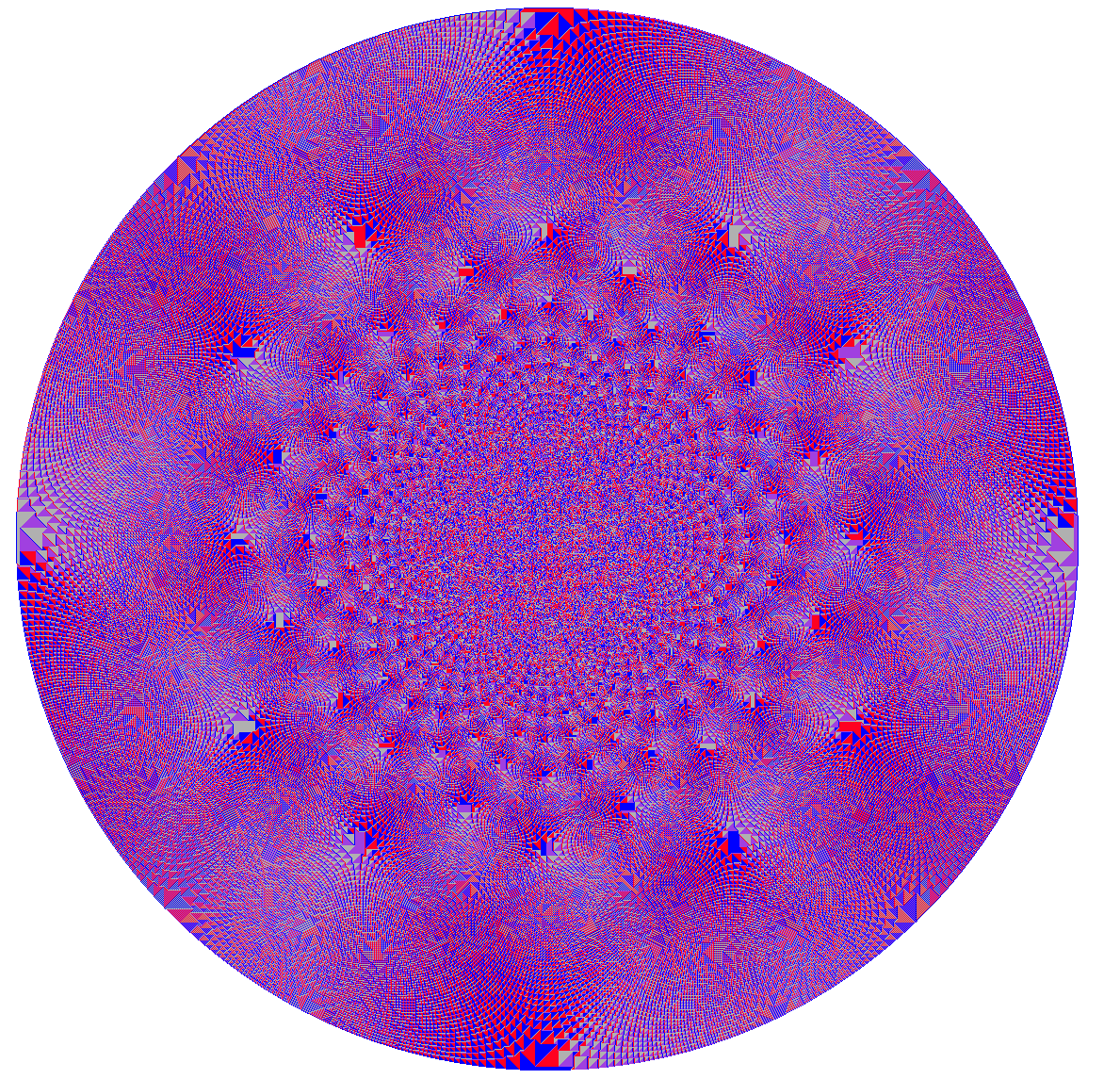}
\caption{Rotor-router aggregate of one million particles in $\Z^2$.  Each site is colored according to the direction of its rotor.}
\label{rotor1m}
\end{figure}

Jim Propp observed from simulations in two dimensions that the regions $A_n$ are extraordinarily close to circular, and asked why this was so \cite{Kleber,Propp}.  Despite the impressive empirical evidence for circularity, the best result known until now \cite{LP} says only that if $A_n$ is rescaled to have unit volume, the volume of the symmetric difference of $A_n$ with a ball of unit volume tends to zero as a power of $n$, as $n \uparrow \infty$.  The main outline of the argument is summarized in \cite{intelligencer}.  Fey and Redig  \cite{FR} also show that $A_n$ contains a diamond.  In particular, these results do not rule out the possibility of ``holes'' in $A_n$ far from the boundary or of long tendrils extending far beyond the boundary of the ball, provided the volume of these features is negligible compared to $n$.

Our main result is the following, which rules out the possibility of holes far from the boundary or of long tendrils in the rotor-router shape.  For $r\geq 0$ let
	\[ B_r = \{ x \in \Z^d ~:~ |x| < r \}. \] 

\begin{theorem}
\label{rotorcirc}
Let $A_n$ be the region formed by rotor-router aggregation in $\Z^d$ starting from $n$ particles at the origin and any initial rotor state.  There exist constants $c,c'$ depending only on $d$, such that
	\[ B_{r-c\log r} \subset A_n \subset B_{r (1+c'r^{-1/d}\log r)} \]
where $r=(n/\omega_d)^{1/d}$, and $\omega_d$ is the volume of the unit ball in $\R^d$.
\end{theorem}

We remark that the same result holds when the rotors evolve according to stacks of bounded discrepancy; see the remark following Lemma~\ref{odomflow}.

Internal diffusion limited aggregation (``internal DLA'') is an analogous growth model defined using random walks instead of rotor-router walks.  Starting with $n$ particles at the origin, each particle in turn performs simple random walk until it reaches an unoccupied site.  Lawler, Bramson and Griffeath \cite{LBG} showed that for internal DLA in $\Z^d$, the occupied region $A_n$, rescaled by a factor of $n^{1/d}$, converges with probability one to a Euclidean ball in $\R^d$ as $n \rightarrow \infty$.  Lawler \cite{Lawler95} estimated the rate of convergence.
By way of comparison with Theorem~\ref{rotorcirc}, if $I_n$ is the internal DLA region formed from $n$ particles started at the origin, the best known bounds \cite{Lawler95} are (up to logarithmic factors)
	\[ B_{r-r^{1/3}} \subset I_n \subset B_{r+r^{1/3}} \]
for all sufficiently large $n$, with probability one.

We also study another model which is slightly more difficult to define, but much easier to analyze.  
In the {\it divisible sandpile}, each site $x \in \Z^d$ starts with a quantity of ``mass'' $\nu_0(x) \in \R_{\geq 0}$.  A site {\it topples} by keeping up to mass~$1$ for itself, and distributing the excess (if any) equally among its neighbors.  Thus if $x$ has mass $m> 1$, then each of the $2d$ neighboring sites gains mass $(m-1)/2d$ when we topple $x$, and $x$ is left with mass $1$; if $m\leq 1$, then no mass moves when we topple $x$.

  Note that individual topplings do not commute; however, the divisible sandpile is ``abelian'' in the following sense.

\begin{prop}
\label{divsandpileabelianproperty}
Let $x_1, x_2, \ldots \in \Z^d$ be a sequence with the property that for any $x \in \Z^d$ there are infinitely many terms $x_k=x$.  Let
	\begin{align*} u_k(x) &= \text{\em total mass emitted by $x$ after toppling $x_1, \ldots, x_k$}; \\
			       \nu_k(x) &= \text{\em amount of mass present at $x$ after toppling $x_1,\ldots, x_k$}.
			       \end{align*}
Then $u_k \uparrow u$ and $\nu_k \to \nu \leq 1$.  Moreover, the limits $u$ and $\nu$ are independent of the sequence $\{x_k\}$.
\end{prop}

\begin{figure}
\centering
\includegraphics[scale=.3333333333]{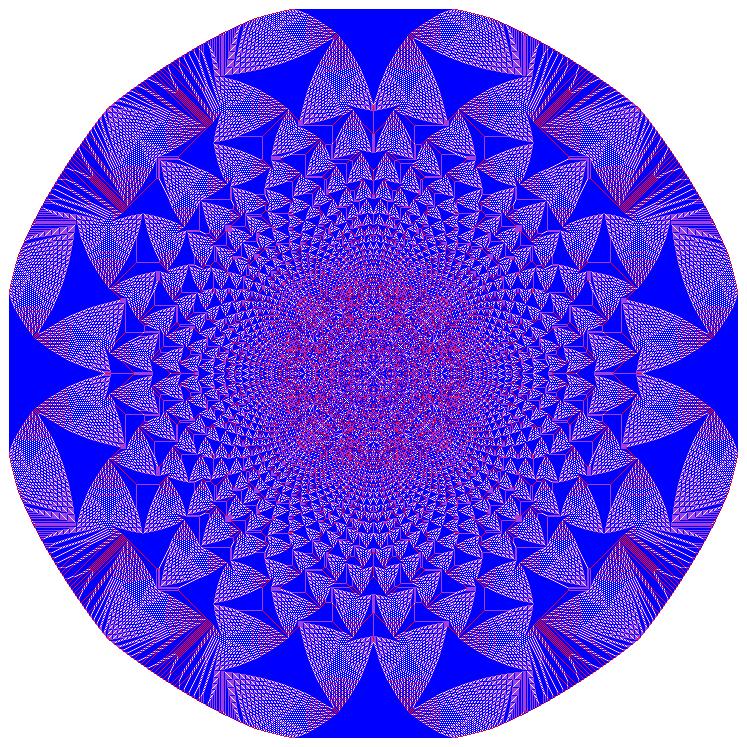}
\caption{Classical abelian sandpile aggregate of one million particles in $\Z^2$.  Colors represent the number of grains at each site.}
\label{sandpile1m}
\end{figure}

The abelian property can be generalized as follows: after performing some topplings, we can add some additional mass and then continue toppling.  The resulting limits $u$ and $\nu$ will be the same as in the case when all mass was initially present.  For a further generalization, see \cite{scalinglimit}.

The limiting function $u$ in Proposition~\ref{divsandpileabelianproperty} is the \emph{odometer function} for the divisible sandpile.  This function plays a central role in our analysis.  The limit $\nu$ represents the final mass distribution.  Sites $x \in \Z^d$ with $\nu(x)=1$ are called {\it fully occupied}.  Proposition~\ref{divsandpileabelianproperty} is proved in section~3, along with the following.

\begin{theorem}
\label{divsandcircintro}
For $m \geq 0$ let $D_m \subset \Z^d$ be the domain of fully occupied sites for the divisible sandpile formed from a pile of mass $m$ at the origin.  There exist constants $c,c'$ depending only on $d$, such that
	\[ B_{r-c} \subset D_m \subset B_{r+c'}, \]
where $r = (m/\omega_d)^{1/d}$ and $\omega_d$ is the volume of the unit ball in $\R^d$.
\end{theorem} 

The divisible sandpile is similar to the ``oil game'' studied by Van den Heuvel \cite{VdH}.  In the terminology of \cite{FR}, it also corresponds to the $h \rightarrow -\infty$ limit of the classical abelian sandpile (defined below), that is, the abelian sandpile started from the initial condition in which every site has a very deep ``hole.''

In the classical abelian sandpile model \cite{BTW}, each site in $\Z^d$ has an integer number of grains of sand; if a site has at least $2d$ grains, it {\it topples}, sending one grain to each neighbor.  If $n$ grains of sand are started at the origin in $\Z^d$, write $S_n$ for the set of sites that are visited during the toppling process; in particular, although a site may be empty in the final state, we include it in $S_n$ if it was occupied at any time during the evolution to the final state. 

\begin{figure}
\centering
\includegraphics[scale=.25]{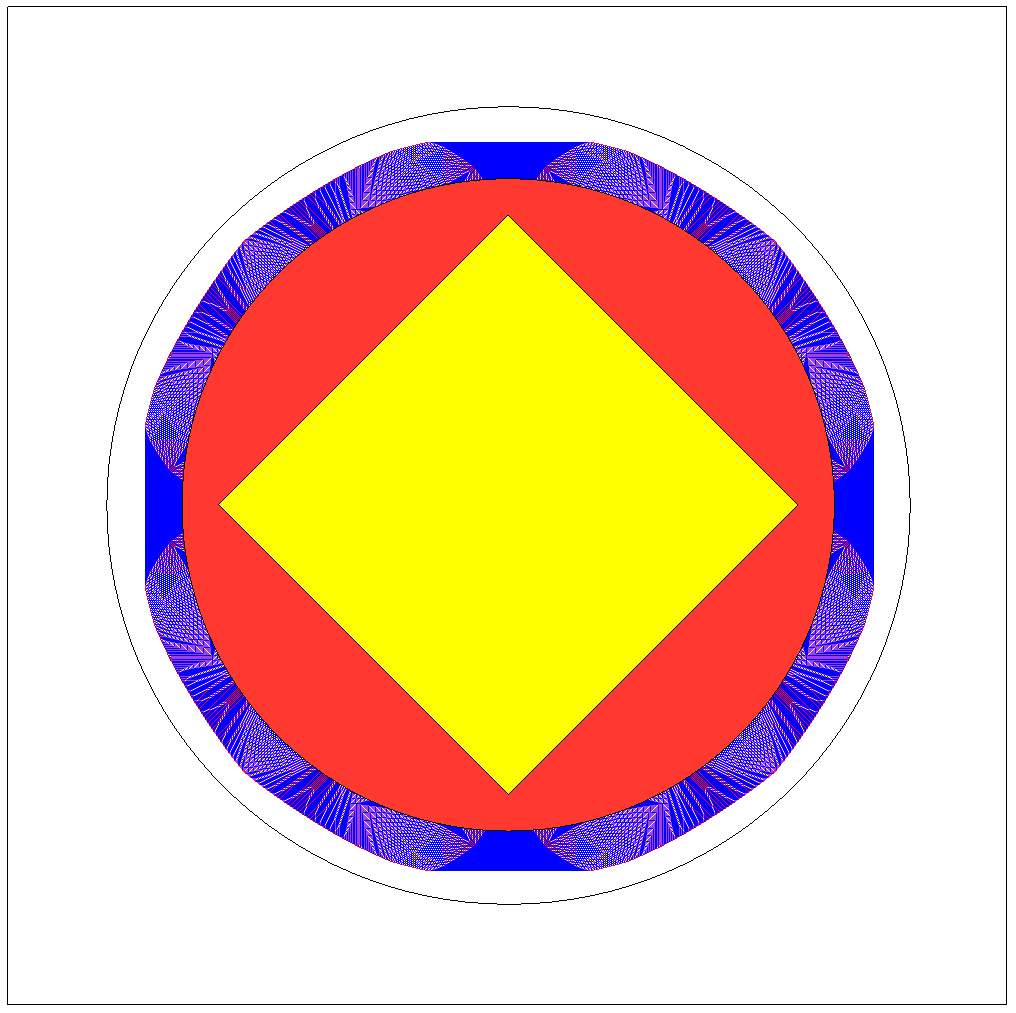}
\caption{Known bounds on the shape of the classical abelian sandpile in $\Z^2$.  The inner diamond and outer square are due to Le Borgne and Rossin \cite{LBR}; the inner and outer circles are those in Theorem~\ref{sandpilecircintro}.}
\end{figure}

Until now the best known constraints on the shape of $S_n$ in two dimensions were due to Le Borgne and Rossin \cite{LBR}, who proved that
	\[ \{x \in \Z^2 \,|\, x_1+x_2 \leq \sqrt{n/12}-1 \} \subset S_n \subset \{x \in \Z^2 \,|\, x_1,x_2 \leq \sqrt{n}/2 \}. \]
Fey and Redig~\cite{FR} proved analogous bounds in higher dimensions, and extended these bounds to arbitrary values of the height parameter $h$.  This parameter is discussed in section~4.

The methods used to prove the near-perfect circularity of the divisible sandpile shape in Theorem~\ref{divsandcircintro} can be used to give constraints on the shape of the classical abelian sandpile, improving on the bounds of~\cite{FR} and~\cite{LBR}.

\begin{theorem}
\label{sandpilecircintro}
Let $S_n$ be the set of sites that are visited by the classical abelian sandpile model in $\Z^d$, starting from $n$ particles at the origin.  Write $n = \omega_d r^d$.  Then for any $\epsilon>0$ we have
	\[ B_{c_1r - c_2} \subset S_n \subset B_{c'_1r + c'_2} \]
where
	\[ c_1 = (2d-1)^{-1/d}, \qquad c'_1 = (d-\epsilon)^{-1/d}. \]
The constant $c_2$ depends only on $d$, while $c'_2$ depends only on $d$ and $\epsilon$.
\end{theorem}

Note that Theorem~\ref{sandpilecircintro} does not settle the question of the asymptotic shape of $S_n$, and indeed it is not clear from simulations whether the asymptotic shape in two dimensions is a disc or perhaps a polygon (Figure~\ref{sandpile1m}).  To our knowledge even the existence of an asymptotic shape is not known.

The rest of the paper is organized as follows.  In section~2, we derive the basic Green's function estimates that are used in the proofs of Theorems~\ref{rotorcirc}, \ref{divsandcircintro} and~\ref{sandpilecircintro}.  In section~3 we prove Proposition~\ref{divsandpileabelianproperty} and Theorem~\ref{divsandcircintro} for the divisible sandpile.  In section~4 we adapt the methods of the previous section to prove Theorem~\ref{sandpilecircintro} for the classical abelian sandpile model.  Section~5 is devoted to the proof of Theorem~\ref{rotorcirc}.

\section{Basic Estimate}

Write $(X_k)_{k\geq 0}$ for simple random walk in $\Z^d$, and for $d\geq 3$ denote by 
	\[ g(x) = \EE_o \# \{ k|X_k=x \} \]
the expected number of visits to $x$ by simple random walk started at the origin.  This is the \emph{discrete harmonic Green's function} in $\Z^d$; it satisfies $\Delta g(x) = 0$ for $x\neq o$, and $\Delta g(o)=-1$, where $\Delta$ is the discrete Laplacian
	\[ \Delta g(x) = \frac{1}{2d} \sum_{y \sim x} g(y) - g(x). \]
The sum is over the $2d$ lattice neighbors $y$ of $x$.  In dimension $d=2$, simple random walk is recurrent, so the expectation defining $g(x)$ is infinite.  Here we define the {\it potential kernel}
	\begin{equation} \label{potentialkerneldef} g(x) = \lim_{n \rightarrow \infty} g_n(x)-g_n(o) \end{equation}
where
	\[ g_n(x) = \EE_o \# \{k \leq n|X_k=x\}. \]
The limit defining $g(x)$ in (\ref{potentialkerneldef}) is finite \cite{Lawler,Spitzer}, and $g(x)$ has Laplacian $\Delta g(x)=0$ for $x\neq o$, and $\Delta g(o)=-1$.  Note that (\ref{potentialkerneldef}) is the negative of the usual definition of the potential kernel; we have chosen this sign convention so that $g$ has the same Laplacian in dimension two as in higher dimensions.

Fix a real number $m>0$ and consider the function on $\Z^d$
	\begin{equation} \label{gammadef} \widetilde{\gamma}_d(x) = |x|^2 + mg(x). \end{equation}
Let $r$ be such that $m = \omega_d r^d$, and let 
	\begin{equation} \gamma_d(x) = \widetilde{\gamma}_d(x)-\widetilde{\gamma}_d(\floor{r}e_1) \label{subtractoffconstant} \end{equation}
where $e_1$ is the first standard basis vector in $\Z^d$.  The function $\gamma_d$ plays a central role in our analysis.  To see where it comes from, recall the divisible sandpile \emph{odometer function} of Proposition~\ref{divsandpileabelianproperty}
	\[ u(x) = \text{total mass emitted from $x$}. \]
Let $D_m \subset \Z^d$ be the domain of fully occupied sites for the divisible sandpile formed from a pile of mass $m$ at the origin.  For $x \in D_m$, since each neighbor $y$ of $x$ emits an equal amount of mass to each of its $2d$ neighbors, we have
	\begin{align*} \Delta u(x) &= \frac{1}{2d} \sum_{y \sim x} u(y) - u(x) \\
					        &= \text{mass received by $x$} - \text{mass emitted by $x$} \\
					        &= 1 - m\delta_{ox}. \end{align*}
Moreover, $u=0$ on $\partial D_m$.  By construction, the function $\gamma_d$ obeys the same Laplacian condition: $\Delta \gamma_d = 1 - m\delta_o$; and as we will see shortly, $\gamma_d \approx 0$ on $\partial B_r$.  Since we expect the domain $D_m$ to be close to the ball $B_r$, we should expect that $u \approx \gamma_d$.  In fact, we will first show that $u$ is close to $\gamma_d$, and then use this to conclude that $D_m$ is close to $B_r$.

We will use the following estimates for the Green's function \cite{FU,Uchiyama}; see also \cite[Theorems 1.5.4 and 1.6.2]{Lawler}.
 \begin{equation}  \label{standardgreensestimate}
  g(x) = \begin{cases}  -\frac{2}{\pi} \log |x| + \kappa + O(|x|^{-2}), & d=2 \\
			 a_d |x|^{2-d} + O(|x|^{-d}), & d\geq 3. \end{cases} \end{equation}
Here $a_d = \frac{2}{(d-2)\omega_d}$, where $\omega_d$ is the volume of the unit ball in $\R^d$, and $\kappa$ is a constant whose value we will not need to know.  For $x \in \Z^d$ we use $|x|$ to denote the Euclidean norm of $x$.  Here and throughout the paper, constants in error terms denoted $O(\cdot)$ depend only on $d$.

We will need an estimate for $\gamma_d$ near the boundary of the ball $B_r$.  We first consider dimension $d=2$.  From (\ref{standardgreensestimate}) we have
	\begin{equation} \label{gammadefn} \widetilde{\gamma}_2(x) = \phi(x) - \kappa m + O(m|x|^{-2}), \end{equation}
where
	\[ \phi(x) = |x|^2 - \frac{2m}{\pi} \log |x|. \]
In the Taylor expansion of $\phi$ about $|x|=r$ 
	\begin{equation} \label{secondordertaylor} \phi(x) = \phi(r) - \phi'(r)(r-|x|) + \frac12 \phi''(t)(r-|x|)^2 \end{equation}
the linear term vanishes, leaving
	\begin{equation} \label{taylordim2} 
	\gamma_2(x) =  \left(1+\frac{m}{\pi t^2} \right)(r-|x|)^2 + O(m|x|^{-2})
	\end{equation}
for some $t$ between $|x|$ and $r$.

In dimensions $d \geq 3$, from (\ref{standardgreensestimate}) we have
	\[ \widetilde{\gamma}_d(x) = |x|^2 + a_d m |x|^{2-d} + O(m|x|^{-d}). \]
Setting $\phi(x) = |x|^2 + a_d m |x|^{2-d}$, the linear term in the Taylor expansion (\ref{secondordertaylor}) of $\phi$ about $|x|=r$ again vanishes, yielding
	\[ \gamma_d(x) = \left( 1 + (d-1)(r/t)^d \right) (r-|x|)^2 + O(m|x|^{-d}) \]
for $t$ between $|x|$ and $r$.  Together with (\ref{taylordim2}), this yields the following estimates in all dimensions $d\geq 2$.

\begin{lemma} \label{gammalowerbound}
Let $\gamma_d$ be given by (\ref{subtractoffconstant}).
For all $x \in \Z^d$ we have
	\begin{equation} \label{quadraticgrowth} \gamma_d(x) \geq (r-|x|)^2 + O\left(\frac{r^d}{|x|^d}\right). \end{equation}
\end{lemma}

\begin{lemma} \label{gammaupperbound}
Let $\gamma_d$ be given by (\ref{subtractoffconstant}).  Then uniformly in $r$,
	\[ \gamma_d(x) = O(1), \qquad x \in B_{r+1}-B_{r-1}. \]	
\end{lemma}

The following lemma is useful for $x$ near the origin, where the error term in (\ref{quadraticgrowth}) blows up.

\begin{lemma} \label{gammanearorigin}
Let $\gamma_d$ be given by (\ref{subtractoffconstant}).  Then for sufficiently large $r$, we have
	\[ \gamma_d(x) > \frac{r^2}{4}, \qquad \forall x \in B_{r/3}. \]
\end{lemma}

\begin{proof}
Since $\gamma_d(x)-|x|^2$ is superharmonic, it attains its minimum in $B_{r/3}$ at a point $z$ on the boundary.  Thus for any $x \in B_{r/3}$
	\[ \gamma_d(x)-|x|^2 \geq \gamma_d(z) - |z|^2, \]
hence by Lemma~\ref{gammalowerbound}
	\[ \gamma_d(x) \geq (2r/3)^2 - (r/3)^2 + O(1) > \frac{r^2}{4}. \qed \]
\renewcommand{\qedsymbol}{}
\end{proof}

Lemmas~\ref{gammalowerbound} and~\ref{gammanearorigin} together imply the following.

\begin{lemma} \label{gammanonnegativeeverywhere}
Let $\gamma_d$ be given by (\ref{subtractoffconstant}).  There is a constant $a$ depending only on $d$, such that $\gamma_d \geq -a$ everywhere.
\end{lemma}

\section{Divisible Sandpile}

Let $\mu_0$ be a nonnegative function on $\Z^d$ with finite support.  
We start with mass $\mu_0(y)$ at each site $y$.  The operation of {\it toppling} a vertex $x$ yields the mass distribution
	\[ T_x \mu_0 = \mu_0 + \alpha(x)\Delta\delta_x \]
where $\alpha(x) = \text{max}(\mu_0(x)-1,0)$ and $\Delta$ is the discrete Laplacian on $\Z^d$.  Thus if $\mu_0(x)\leq 1$ then $T_x \mu_0 = \mu_0$ and no mass topples; if $\mu_0(x)>1$ then the mass in excess of $1$ is distributed equally among the neighbors of $x$.

Let $x_1, x_2, \ldots \in \Z^d$ be a sequence with the property that for any $x \in \Z^d$ there are infinitely many terms $x_k=x$.  Let
 	\[ \mu_k(y) = T_{x_k} \ldots T_{x_1} \mu_0(y). \]
be the amount of mass present at $y$ after toppling the sites $x_1, \ldots, x_k$ in succession.  The total mass emitted from $y$ during this process is
	\begin{equation} \label{finitetimeodom} u_k(y) := \sum_{j\leq k : x_j=y} \mu_{j-1}(y)-\mu_j(y) = \sum_{j\leq k : x_j=y} \alpha_j(y)
	\end{equation}
where $\alpha_j(y) = \text{max}(\mu_j(y)-1,0)$.

\begin{lemma}
\label{existenceoflimits}
As $k\uparrow \infty$ the functions $u_k$ and $\mu_k$ tend to limits $u_k \uparrow u$ and $\mu_k \rightarrow \mu$.  Moreover, these limits satisfy
	\[ \mu = \mu_0 + \Delta u \leq 1. \]
\end{lemma}
	
\begin{proof}
Write $M = \sum_x \mu_0(x)$ for the total starting mass, and let $B \subset \Z^d$ be a ball centered at the origin containing all points within $L^1$-distance $M$ of the support of $\mu_0$.  Note that if $\mu_k(x)>0$ and $\mu_0(x)=0$, then $x$ must have received its mass from a neighbor, so $\mu_k(y)\geq 1$ for some $y\sim x$.  Since $\sum_x \mu_k(x) = M$, it follows that $\mu_k$ is supported on $B$.  Let $R$ be the radius of $B$, and consider the {\it quadratic weight}
	\[ Q_k = \sum_{x \in \Z^d} \mu_k(x) |x|^2 \leq MR^2. \]
Since $\mu_k(x_k)-\mu_{k-1}(x_k) = -\alpha_k(x_k)$ and for $y\sim x_k$ we have $\mu_k(y)-\mu_{k-1}(y) = \frac{1}{2d} \alpha_k(x_k)$, we obtain
	\[ Q_k - Q_{k-1} = \alpha_k(x_k) \left( \frac{1}{2d} \sum_{y \sim x_k} |y|^2 - |x_k|^2 \right) = \alpha_k(x_k). \]
Summing over $k$ we obtain from (\ref{finitetimeodom})
	\[ Q_k = Q_0 + \sum_{x \in \Z^d} u_k(x). \]
Fixing $x$, the sequence $u_k(x)$ is thus increasing and bounded above, hence convergent.

Given neighboring vertices $x\sim y$, since $y$ emits an equal amount of mass to each of its $2d$ neighbors, it emits mass $u_k(y)/2d$ to $x$ up to time $k$.  Thus $x$ receives a total mass of $\frac{1}{2d} \sum_{y\sim x} u_k(y)$ from its neighbors up to time $k$.  Comparing the amount of mass present at $x$ before and after toppling, we obtain
	\[ \mu_k(x) = \mu_0(x) + \Delta u_k(x). \]
Since $u_k \uparrow u$ we infer that $\mu_k \rightarrow \mu:= \mu_0 + \Delta u$.  Note that if $x_k=x$, then $\mu_k(x) \leq 1$.  Since for each $x \in \Z^d$ this holds for infinitely many values of $k$, the limit satisfies $\mu \leq 1$.
\end{proof}

A function $s$ on $\Z^d$ is {\it superharmonic} if $\Delta s \leq 0$.  Given a function $\gamma$ on $\Z^d$ the {\it least superharmonic majorant} of $\gamma$ is the function
	\[ s(x) = \text{inf} \{ f(x)~|~ f \text{ is superharmonic and } f\geq \gamma \}. \]
The study of the least superharmonic majorant is a classical topic in analysis and PDE; see, for example, \cite{Koosis}.  Note that if $f$ is superharmonic and $f \geq \gamma$ then
	\[ f(x) \geq \frac{1}{2d} \sum_{y \sim x} f(y) \geq \frac{1}{2d} \sum_{y \sim x} s(y). \]
Taking the infimum on the left side we obtain that $s$ is superharmonic.

\begin{lemma}
\label{discretemajorant}
The limit $u$ in Lemma~\ref{existenceoflimits} is given by $u = s + \gamma$, where
	\[ \gamma(x) = |x|^2 + \sum_{y \in \Z^d} g(x-y) \mu_0(y) \]
and $s$ is the least superharmonic majorant of $-\gamma$.  
\end{lemma}

\begin{proof}
By Lemma~\ref{existenceoflimits} we have
	\[ \Delta u = \mu - \mu_0 \leq 1-\mu_0. \]
Since $\Delta \gamma = 1-\mu_0$, the difference $u-\gamma$ is superharmonic.  As $u$ is nonnegative, it follows that $u-\gamma \geq s$.  For the reverse inequality, note that $s+\gamma-u$ is superharmonic on the domain $D = \{x~|~\mu(x)=1\}$ of fully occupied sites and is nonnegative outside $D$, hence nonnegative inside $D$ as well. 
\end{proof}

As a corollary of Lemmas~\ref{existenceoflimits} and~\ref{discretemajorant}, we obtain the abelian property of the divisible sandpile, Proposition~\ref{divsandpileabelianproperty}, which was stated in the introduction.

We now turn to the case of a point source mass $m$ started at the origin: $\mu_0=m\delta_{o}$.  More general starting distributions are treated in \cite{scalinglimit}, where we identify the scaling limit of the divisible sandpile model and show that it coincides with that of internal DLA and of the rotor-router model.  In the case of a point source of mass $m$, the natural question is to identify the shape of the resulting domain $D_m$ of fully occupied sites, i.e.\ sites $x$ for which $\mu(x)=1$.  According to Theorem~\ref{divsandcirc}, $D_m$ is extremely close to a ball of volume $m$; in fact, the error in the radius is a constant independent of $m$.
As before, for $r\geq 0$ we write
	\[ B_r = \{ x \in \Z^d ~:~ |x| < r \} \]
for the lattice ball of radius $r$ centered at the origin.
\begin{theorem}
\label{divsandcirc}
For $m \geq 0$ let $D_m \subset \Z^d$ be the domain of fully occupied sites for the divisible sandpile formed from a pile of size $m$ at the origin.  There exist constants $c,c'$ depending only on $d$, such that
	\[ B_{r-c} \subset D_m \subset B_{r+c'}, \]
where $r = (m/\omega_d)^{1/d}$ and $\omega_d$ is the volume of the unit ball in $\R^d$.
\end{theorem} 


The idea of the proof is to use Lemma~\ref{discretemajorant} along with the basic estimates on $\gamma$, Lemmas~\ref{gammalowerbound} and~\ref{gammaupperbound}, to obtain estimates on the odometer function
	\[ u(x) = \text{total mass emitted from $x$}. \]
We will need the following simple observation.
	
\begin{lemma}
\label{boundarypath}
For every point $x\in D_m - \{o\}$ there is a path $x=x_0 \sim x_1 \sim \ldots \sim x_k=o$ in $D_m$ with $u(x_{i+1}) \geq u(x_i)+1$.
\end{lemma}
\begin{proof}
If $x_i \in D_m - \{o\}$, let $x_{i+1}$ be a neighbor of $x_i$ maximizing $u(x_{i+1})$.  Then $x_{i+1} \in D_m$ and
	\begin{align*} u(x_{i+1}) &\geq \frac{1}{2d} \sum_{y \sim x_i} u(y) \\
		&= u(x_i) + \Delta u(x_i) \\
		&= u(x_i)+1,  \end{align*}
where in the last step we have used the fact that $x_i \in D_m$.
\end{proof}
   
\begin{proof}[Proof of Theorem~\ref{divsandcirc}]
We first treat the inner estimate.  Let $\gamma_d$ be given by (\ref{subtractoffconstant}).
By Lemma~\ref{discretemajorant} the function $u-\gamma_d$ is superharmonic,  so its minimum in the ball $B_r$ is attained on the boundary.  Since $u\geq 0$, we have by  Lemma~\ref{gammaupperbound}
	\[ u(x)-\gamma_d(x) \geq -C, \qquad x \in \partial B_r \]
for a constant $C$ depending only on $d$.  Hence by Lemma~\ref{gammalowerbound},
	\begin{equation} \label{odomlowerbound} u(x) \geq  (r-|x|)^2 - C'r^d/|x|^d, \qquad x \in B_r. \end{equation}
for a constant $C'$ depending only on $d$.
It follows that there is a constant $c$, depending only on $d$, such that $u(x)>0$ whenever $r/3\leq |x|<r-c$.  Thus $B_{r-c}-B_{r/3} \subset D_m$.  For $x \in B_{r/3}$, by Lemma~\ref{gammanearorigin} we have $u(x) \geq r^2/4-C>0$, hence $B_{r/3} \subset D_m$.

For the outer estimate, note that $u-\gamma_d$ is harmonic on $D_m$.  By Lemma~\ref{gammanonnegativeeverywhere} we have $\gamma_d \geq -a$ everywhere, where $a$ depends only on $d$.  Since $u$ vanishes on $\partial D_m$ it follows that $u-\gamma_d \leq a$ on $D_m$.  
Now for any $x \in D_m$ with $r-1<|x|\leq r$, we have by Lemma~\ref{gammaupperbound}
	\[ u(x) \leq \gamma_d(x) + a \leq c' \]
for a constant $c'$ depending only on $d$.  Lemma~\ref{boundarypath} now implies that $D_m \subset B_{r+c'+1}$.
\end{proof}

\section{Classical Sandpile}

We consider a generalization of the classical abelian sandpile, proposed by Fey and Redig \cite{FR}.  Each site in $\Z^d$ begins with a ``hole'' of depth $H$.  Thus, each site absorbs the first $H$ grains it receives, and thereafter functions normally, toppling once for each additional $2d$ grains it receives.  If $H$ is negative, we can interpret this as saying that every site starts with $h=-H$ grains of sand already present.  Aggregation is only well-defined in the regime $h\leq 2d-2$, since for $h=2d-1$ the addition of a single grain already causes every site in $\Z^d$ to topple infinitely often.

Let $S_{n,H}$ be the set of sites that are visited if $n$ particles start at the origin in $\Z^d$.  Fey and Redig \cite[Theorem~4.7]{FR} prove that
	\[ \lim_{H \to \infty} \limsup_{n \to \infty}  \frac{H}{n} \# \left( S_{n,H} \bigtriangleup B_{H^{-1/d}r} \right) = 0, \]
where $n = \omega_d r^d$, and $\bigtriangleup$ denotes symmetric difference.  The following theorem strengthens this result.

\begin{theorem}
\label{sandpilecirc}
Fix an integer $H \geq 2-2d$.  Let $S_n = S_{n,H}$ be the set of sites that are visited by the classical abelian sandpile model in $\Z^d$, starting from $n$ particles at the origin, if every lattice site begins with a hole of depth $H$.  Write $n = \omega_d r^d$.  Then
	\[ B_{c_1r - c_2} \subset S_{n,H} \]
where
	\[ c_1 = (2d-1+H)^{-1/d} \]
and $c_2$ is a constant depending only on $d$.  Moreover if $H \geq 1-d$, then for any $\epsilon>0$ we have
	 \[ S_{n,H} \subset B_{c'_1r + c'_2} \]
where
 	\[ c'_1 = (d-\epsilon+H)^{-1/d} \]
and $c'_2$ is independent of $n$ but may depend on $d$, $H$ and $\epsilon$.
\end{theorem}

Note that the ratio $c_1/c'_1 \uparrow 1$ as $H \uparrow \infty$.  Thus, the classical abelian sandpile run from an initial state in which each lattice site starts with a deep hole yields a shape very close to a ball.  Intuitively, one can think of the classical sandpile with deep holes as approximating the divisible sandpile, whose limiting shape is a ball by Theorem~\ref{divsandcirc}.  Following this intuition, we can adapt the proof of Theorem~\ref{divsandcirc} to prove Theorem~\ref{sandpilecirc}; just one additional averaging trick is needed, which we explain below.

Consider the odometer function for the abelian sandpile
	\[  u(x) = \text{total number of grains emitted from }x. \]
Let $T_n = \{x|u(x)>0\}$ be the set of sites which topple at least once.  Then 
	\[ T_n \subset S_n \subset  T_n \cup \partial T_n. \]
In the final state, each site which has toppled retains between $0$ and $2d-1$ grains, in addition to the $H$ that it absorbed.  Hence
	\begin{equation} \label{dumblaplbounds} H \leq \Delta u(x) + n\delta_{ox} \leq 2d-1+H, \qquad x \in T_n. \ 
	\end{equation}
We can improve the lower bound by averaging over a small box.  For $x \in \Z^d$ let
	\[ Q_k(x) = \{ y\in \Z^d \,:\, ||x-y||_{\infty} \leq k \} \]
be the box of side length $2k+1$ centered at $x$, and let
	\[ u^{(k)}(x) = (2k+1)^{-d} \sum_{y \in Q_k(x)} u(y). \]
Write
	\[ T_n^{(k)} = \{ x \,|\, Q_k(x) \subset T_n \}. \]
Le Borgne and Rossin \cite{LBR} observe that if $T$ is a set of sites all of which topple, the number of grains remaining in $T$ is at least the number of edges internal to $T$: indeed, for each internal edge, the endpoint that topples last sends the other a grain which never moves again.  Since the box $Q_k(x)$ has $2dk(2k+1)^{d-1}$ internal edges, we have
	\begin{equation} \label{averagedinabox} \Delta u^{(k)}(x) \geq \frac{2k}{2k+1}d +H - \frac{n}{(2k+1)^d} \textbf{1}_{Q_k(o)}(x), \qquad x \in T_n^{(k)}. \end{equation}
The following lemma is analogous to Lemma~\ref{boundarypath}.
 
\begin{lemma}
\label{boundarypathology}
For every point $x\in T_n$ adjacent to $\partial T_n$ there is a path $x=x_0 \sim x_1 \sim \ldots \sim x_m = o$ in $T_n$ with $u(x_{i+1}) \geq u(x_i)+1$.
\end{lemma}
 
\begin{proof}
By (\ref{dumblaplbounds}) we have
	\[ \frac{1}{2d} \sum_{y \sim x_i} u(y) \geq u(x_i). \]
Since $u(x_{i-1})<u(x_i)$, some term $u(y)$ in the sum above must exceed $u(x_i)$.  Let $x_{i+1}=y$.
\end{proof} 
 
\begin{proof}[Proof of Theorem~\ref{sandpilecirc}]
Let
	\[ \widetilde{\xi}_d(x) = (2d-1+H)|x|^2 + ng(x), \]
and let 
	\[ \xi_d(x) = \widetilde{\xi}_d(x) - \widetilde{\xi}_d(\floor{c_1r}e_1). \]
Taking $m = n/(2d-1+H)$ in Lemma~\ref{gammaupperbound}, we have
	\begin{equation} \label{minattainedonboundary} u(x)-\xi_d(x) \geq -\xi_d(x) \geq -C(2d-1+H), \qquad x\in\partial B_{c_1r} \end{equation}
for a constant $C$ depending only on $d$.
By (\ref{dumblaplbounds}), $u-\xi_d$ is superharmonic, so $u-\xi_d \geq -C(2d-1+H)$ in all of $B_{c_1r}$.  Hence by Lemma~\ref{gammalowerbound} we have for $x \in B_{c_1 r}$
	\begin{equation} \label{odomlowerboundagain} u(x) \geq (2d-1+H) \left( (c_1r - |x|)^2 - C'(c_1 r)^d/|x|^d \right), 
	\end{equation}
where $C'$ depends only on $d$.
It follows that $u$ is positive on $B_{c_1r - c_2}-B_{c_1 r/3}$ for a suitable constant $c_2$ depending only on $d$.  For $x \in B_{c_1 r/3}$, by Lemma~\ref{gammanearorigin} we have $u(x) > (2d-1+H)(c_1^2 r^2 / 4 -C) > 0$.  Thus $B_{c_1 r - c_2} \subset T_n \subset S_n$.  

For the outer estimate, let
	\[ \hat{\psi}_d(x) = (d-\epsilon+H)|x|^2 + ng(x). \]
Choose $k$ large enough so that $\frac{2k}{2k+1} d \geq d-\epsilon$, and define
	\[ \widetilde{\psi}_d(x) = (2k+1)^{-d} \sum_{y \in Q_k(x)} \hat{\psi}_d(y). \]
Finally, let
	\[ \psi_d(x) = \widetilde{\psi}_d(x) - \widetilde{\psi}_d(\floor{c'_1r}e_1). \]
By (\ref{averagedinabox}), $u^{(k)}-\psi_d$ is subharmonic on $T_n^{(k)}$.  Taking $m=n/(d-\epsilon+H)$ in Lemma~\ref{gammanonnegativeeverywhere}, there is a constant $a$ depending only on $d$, such that $\psi_d \geq -a(d+H)$ everywhere.  Since $u^{(k)} \leq (2d+H)^{(d+1)k}$ on $\partial T_n^{(k)}$ it follows that $u^{(k)}-\psi_d \leq a(d+H)+(2d+H)^{(d+1)k}$ on $T_n^{(k)}$.  Now for any $x \in S_n$ with $c'_1r-1 < |x| \leq c'_1 r$ we have by Lemma~\ref{gammaupperbound}
	\[ u^{(k)}(x) \leq \psi_d(x) + a(d+H) + (2d+H)^{(d+1)k} \leq \tilde{c}_2 \]
for a constant $\tilde{c}_2$ depending only on $d$, $H$ and $\epsilon$.  Then $u(x) \leq c'_2 := (2k+1)^d \tilde{c}_2$.  Lemma~\ref{boundarypathology} now implies that $T_n \subset B_{c'_1 r + c'_2}$, and hence 
	\[ S_n \subset T_n \cup \partial T_n \subset B_{c'_1 r + c'_2+1}. \qed \]
\renewcommand{\qedsymbol}{}
\end{proof}

We remark that the crude bound of $(2d+H)^{(d+1)k}$ used in the proof of the outer estimate can be improved to a bound of order $k^2H$, and the final factor of $(2k+1)^d$ can be replaced by a constant factor independent of $k$ and $H$, using the fact that a nonnegative function on $\Z^d$ with bounded Laplacian cannot grow faster than quadratically; see \cite{scalinglimit}.

\section{Rotor-Router Model}

Given a function $f$ on $\Z^d$, for a directed edge $(x,y)$ write
	\[ \nabla f (x,y) = f(y)-f(x). \]
Given a function $s$ on directed edges in $\Z^d$, write
	\[ \div s (x) = \frac{1}{2d} \sum_{y \sim x} s(x,y). \]
The discrete Laplacian of $f$ is then given by
	\[ \Delta f(x) = \div \nabla f = \frac{1}{2d} \sum_{y \sim x} f(y)-f(x). \]

\subsection{Inner Estimate}

Fixing $n \geq 1$, consider the odometer function for rotor-router aggregation
	\[ u(x) = \text{total number exits from $x$ by the first $n$ particles}. \]
We learned the idea of using the odometer function to study the rotor-router shape from Matt Cook \cite{Cook}.

\begin{lemma}
\label{odomflow}
For a directed edge $(x,y)$ in $\Z^d$, denote by $\kappa(x,y)$ the net number of crossings 
from $x$ to $y$ performed by the first $n$ particles in rotor-router aggregation.  
Then
	\begin{equation} \label{gradodom} \nabla u(x,y) = -2d\kappa(x,y) + R(x,y) 
\end{equation}
for some edge function $R$ which satisfies
	\[ |R(x,y)| \leq 4d-2 \]
for all edges $(x,y)$.
\end{lemma}

\begin{remark}
In the more general setting of rotor stacks of bounded discrepancy, the $4d-2$ will be replaced by a different constant here.
\end{remark}

\begin{proof}
Writing $N(x,y)$ for the number of particles routed from $x$ to $y$, we have
	\[ \frac{u(x)-2d+1}{2d} \leq N(x,y) \leq \frac{u(x)+2d-1}{2d} \]
hence
	\begin{align*}  | \nabla u(x,y) + 2d \kappa(x,y) | &= | u(y) - u(x) + 2dN(x,y) - 2dN(y,x) | \\
		&\leq 4d - 2. \qed \end{align*} 
\renewcommand{\qedsymbol}{} 
\end{proof}

In what follows, $C_0, C_1, \ldots$ denote constants depending only on $d$.

\begin{lemma}
\label{diameterbound}
Let $\Omega \subset \Z^d-\{o\}$ with $2 \leq \#\Omega < \infty$.  Then
	\[ \sum_{y \in \Omega} |y|^{1-d} \leq C_0\, \text{\em Diam}(\Omega). \]
\end{lemma}

\begin{proof}
For each positive integer $k$, let
	\[ \mathcal{S}_k = \{y \in \Z^d \,:\, k\leq |y| < k+1 \}. \]
Then 
	\[ \sum_{y \in S_k} |y|^{1-d} \leq k^{1-d} \# \mathcal{S}_k \leq C'_0 \]
for a constant $C'_0$ depending only on $d$.
Since $\Omega$ can intersect at most $\Diam (\Omega)+1 \leq 2 \Diam(\Omega)$ distinct sets $\mathcal{S}_k$, taking $C_0 = 2C'_0$ the proof is complete.
\end{proof}

\begin{lemma}
\label{letsgetthispartystarted}
Let $G = G_{B_r}$ be the Green's function for simple random walk in $\Z^d$ stopped on exiting $B_r$.  For any $\rho\geq 1$ and $x \in B_r$,
	\begin{equation} \label{theparty}
	\sum_{\substack{y \in B_r \\ |x-y| \leq \rho}} \sum_{z \sim  y} |G(x,y)-G(x,z)| \leq C_1\rho.
	\end{equation}
\end{lemma}

\begin{proof}
Let $(X_t)_{t\geq 0}$ denote simple random walk in $\Z^d$,  and let $T$ be the first exit time from $B_r$. For fixed $y$, the function 
	\begin{equation} \label{stoppedgreen} A(x) = g(x-y) - \EE_x g(X_T-y) \end{equation}
has Laplacian $\Delta A(x) = - \delta_{xy}$ in $B_r$ and vanishes on $\partial B_r$, hence $A(x)=G(x,y)$.  

Let $x,y \in B_r$ and $z\sim y$.  From (\ref{standardgreensestimate}) we have
	\[ |g(x-y)-g(x-z)| \leq\frac{C_2}{|x-y|^{d-1}},  \qquad  y,z \neq x. \]  
Using the triangle inequality together with (\ref{stoppedgreen}), we obtain
	\begin{align*} |G(x,y)-G(x,z)| &\leq |g(x-y)-g(x-z)| + \EE_x | g(X_T-y)-g(X_T-z)| \\	
		&\leq  \frac{C_2}{|x-y|^{d-1}} + \sum_{w\in \partial B_r} H_x(w) \frac{C_2}{|w-y|^{d-1}}, \end{align*}
where $H_x(w) = \PP_x(X_T = w)$.

Write $D = \{y \in B_r \,:\, |x-y|\leq \}$.  Then
	\begin{equation} \label{twointegrals}  \sum_{\substack{y \in D \\ y \neq x}} \sum_{\substack{z \sim  y \\ z\neq x}} |G(x,y)-G(x,z)| \leq C_3 \rho + C_2 \sum_{w \in \partial B_r} H_x(w) \sum_{y \in D} |w-y|^{1-d}.
	\end{equation}
Taking $\Omega = w-D$ in Lemma~\ref{diameterbound}, the inner sum on the right is at most $C_0 $Diam$(D) \leq 2C_0\rho$, so the right side of (\ref{twointegrals}) is bounded above by $C_1 \rho$ for a suitable $C_1$.

Finally, the terms in which $y$ or $z$ coincides with $x$ make a negligible contribution to the sum in (\ref{theparty}), since for $y\sim x \in \Z^d$
	\[ |G(x,x)-G(x,y)| \leq |g(o)-g(x-y)| + \EE_x |g(X_T-x)-g(X_T-y)| \leq C_4. \qed \]
\renewcommand{\qedsymbol}{}
\end{proof}

\begin{figure}
\label{hitthehyperplanefigure}
\centering
\resizebox{3in}{!}{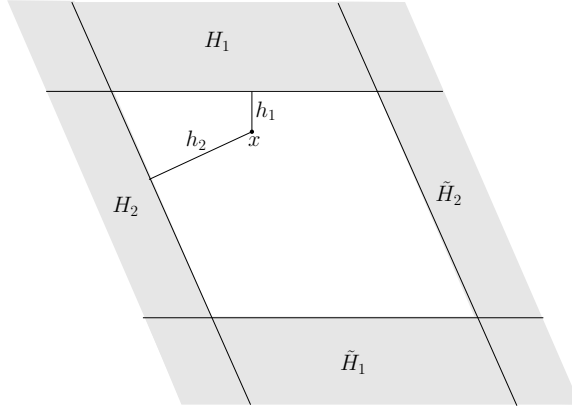}
\caption{Diagram for the Proof of Lemma~\ref{hitthehyperplane}.}
\end{figure}


\begin{lemma}
\label{hitthehyperplane}
Let $H_1, H_2$ be linear half-spaces in $\Z^d$, not necessarily parallel to the coordinate axes.  Let $T_i$ be the first hitting time of $H_i$.  If $x \notin H_1 \cup H_2$, then
	\[ \PP_x(T_1>T_2) \leq \frac52 \frac{h_1+1}{h_2} \left(1+\frac{1}{2h_2}\right)^2 \]
where $h_i$ is the distance from $x$ to $H_i$.
\end{lemma}
 
\begin{proof}
If one of $H_1, H_2$ contains the other, the result is vacuous.  Otherwise, let $\widetilde{H}_i$ be the half-space shifted parallel to $H_i^c$ by distance $2h_2$ in the direction of $x$, and let $\widetilde{T}_i$ be the first hitting time of $H_i \cup \widetilde{H}_i$.  Let $(X_t)_{t\geq 0}$ denote simple random walk in $\Z^d$, and write $M_t$ for the (signed) distance from $X_t$ to the hyperplane defining the boundary of $H_1$, with $M_0=h_1$.   Then $M_t$ is a martingale with bounded increments.
Since $\EE_x \widetilde{T}_1 < \infty$, we obtain from optional stopping
	\[ h_1 = \EE_x M_{\widetilde{T}_1} \geq 2h_2 \PP_x(X_{\widetilde{T}_1} \in \widetilde{H}_1) - \PP_x(X_{\widetilde{T}_1} \in H_1), \]
hence
	\begin{equation} \label{justgamblersruin} \PP_x(X_{\widetilde{T}_1} \in \widetilde{H}_1) \leq \frac{h_1+1}{2h_2}. \end{equation}
Likewise, $dM_t^2 - t$ is a martingale with bounded increments, giving
	\begin{align} \label{timeupperbound} \EE_x \widetilde{T}_1 &\leq d\EE_x M_{\widetilde{T}_1}^2 \nonumber \\
		&\leq d (2h_2+1)^2 \PP_x(X_{\widetilde{T}_1} \in \widetilde{H}_1) \nonumber \\
		&\leq d (h_1+1) (2h_2+1)\left(1+\frac{1}{2h_2}\right). \end{align}
	
Let $T = \text{min}(\widetilde{T}_1, \widetilde{T}_2)$.  Denoting by $D_t$ the distance from $X_t$ to the hyperplane defining the boundary of $H_2$, the quantity 
	\[ N_t = \frac{d}{2} \big( D_t^2 + (2h_2 - D_t)^2 \big) - t \]
is a martingale.  Writing $p = \PP_x(T=\widetilde{T}_2)$ we have
	\begin{align*} dh_2^2 = \EE N_0 = \EE N_T &\geq p \frac{d}{2}(2h_2)^2  + (1-p)dh_2^2 - \EE_x T \\
			&\geq (1+p)dh_2^2 - \EE_x T
	\end{align*}  
hence by (\ref{timeupperbound})
	\[ p \leq \frac{\EE_x T}{dh_2^2} \leq 2 \frac{h_1+1}{h_2} \left(1+\frac{1}{2h_2}\right)^2. \]
Finally by (\ref{justgamblersruin})
	\[ \PP(T_1>T_2) \leq p + \PP(X_{\widetilde{T}_1} \in \widetilde{H}_1) \leq \frac52 \frac{h_1+1}{h_2} \left(1+\frac{1}{2h_2}\right)^2. \qed \]
\renewcommand{\qedsymbol}{}
\end{proof} 


\begin{figure}
\centering
\label{geoshellsfig}
\resizebox{3in}{!}{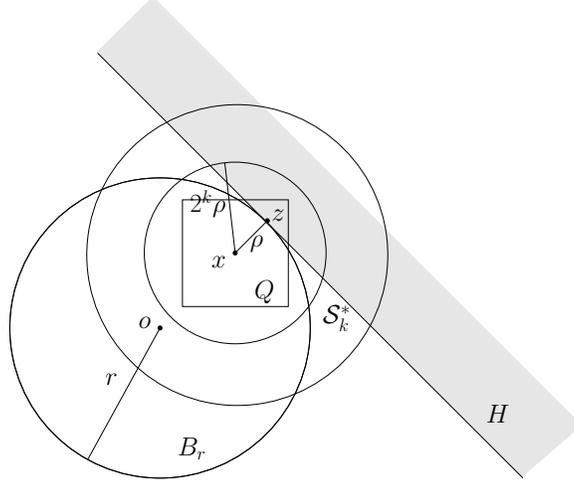}
\caption{Diagram for the proof of Lemma~\ref{geoshells}.}
\end{figure}

\begin{lemma}
\label{geoshells}
Let $x \in B_r$ and let $\rho = r+1-|x|$.  Let
	\begin{equation} \label{shells} \mathcal{S}^*_k = \{ y \in B_r ~:~ 2^k \rho < |x-y| \leq 2^{k+1} \rho \}. \end{equation}
Let $\tau_k$ be the first hitting time of $\mathcal{S}^*_k$, and $T$ the first exit time from $B_r$.  Then 
	\[ \PP_x(\tau_k < T) \leq C_2 2^{-k}. \]
\end{lemma}

\begin{proof}
Let $H$ be the outer half-space tangent to $B_r$ at the point $z\in \partial B_r$ closest to $x$.  Let $Q$ be the cube of side length $2^k \rho/\sqrt{d}$ centered at $x$.  Then $Q$ is disjoint from $\mathcal{S}^*_k$, hence
	\[ \PP_x(\tau_k < T) \leq \PP_x(T_{\partial Q}<T) \leq \PP_x(T_{\partial Q} < T_H) \]
where $T_{\partial Q}$ and $T_H$ are the first hitting times of $\partial Q$ and $H$.  Let $H_1, \ldots, H_{2d}$ be the outer half-spaces defining the faces of $Q$, so that $Q = H_1^c \cap \ldots \cap H_{2d}^c$.  By Lemma~\ref{hitthehyperplane} we have
	\begin{align*}  \PP_x(T_{\partial Q} < T_H) &\leq \sum_{i=1}^{2d} \PP_x(T_{H_i} < T_H) \\
			&\leq \frac52 \sum_{i=1}^{2d} \frac{\text{dist}(x,H)+1}{\text{dist}(x,H_i)} \left(1+\frac{1}{2\,\text{dist}(x,H_i)}\right)^2. \end{align*}
Since dist$(x,H) = |x-z| \leq \rho$ and dist$(x,H_i) = 2^{k-1} \rho/\sqrt{d}$, and $\rho \geq 1$, taking $C_2 = 20\,d^{3/2}(1+\sqrt{d})^2$ completes the proof.
\end{proof}


\begin{lemma}
\label{greengradbound}
Let $G=G_{B_r}$ be the Green's function for random walk stopped on exiting $B_r$.  Let $x \in B_r$ and let $\rho = r+1-|x|$.   Then
 	\[ \sum_{y\in B_r} \sum_{z\sim y} |G(x,y)-G(x,z)| \leq C_3 \rho \log \frac{r}{\rho}. \]
\end{lemma}

\begin{proof}
Let $\mathcal{S}^*_k$ be given by (\ref{shells}), and let 
	\[ W = \{ w \in \partial( \mathcal{S}^*_k \cup \partial \mathcal{S}^*_k) ~:~ |w-x|<2^k \rho \} \]
be the portion of the boundary of the enlarged spherical shell $\mathcal{S}^*_k \cup \partial \mathcal{S}^*_k$ lying closer to $x$.
Let $\tau_W$ be the first hitting time of $W$, and $T$ the first exit time from $B_r$.  For $w\in W$ let 
	\[ H_x(w) = \PP_x(X_{\tau_W \wedge T} = w). \]
For any $y \in \mathcal{S}^*_k$ and $z \sim y$, simple random walk started at $x$ must hit $W$ before hitting either $y$ or $z$, hence
	\[ |G(x,y) - G(x,z)| \leq \sum_{w \in W}  H_x(w)|G(w,y)-G(w,z)|. \]
For any $y \in \mathcal{S}^*_k$ and any $w \in W$ we have 
	\[ |y-w| \leq |y-x| + |w-x| \leq 3\cdot 2^k \rho. \]
Lemma~\ref{letsgetthispartystarted} yields
	\[ \sum_{y \in \mathcal{S}^*_k} \sum_{z \sim y} |G(x,y) - G(x,z)| \leq 3 C_1 2^k \rho \sum_{w \in W} H_x(w). \]
By Lemma~\ref{geoshells} we have $\sum_{w \in W} H_x(w) \leq C_2 2^{-k}$, so the above sum is at most  $3 C_1 C_2 \rho$.  Since the union of shells $\mathcal{S}^*_0, \mathcal{S}^*_1, \ldots, \mathcal{S}^*_{\ceil{\log_2(r/\rho)}}$ covers all of $B_r$ except for those points $y$ within distance $\rho$ of $x$, and $\sum_{|y-x|\leq\rho} \sum_{z \sim y} |G(x,y)-G(x,z)| \leq C_1 \rho$  by Lemma~\ref{letsgetthispartystarted}, the result follows.
\end{proof}

\begin{proof}[Proof of Theorem~\ref{rotorcirc}, Inner Estimate]
Let $\kappa$ and $R$ be defined as in Lemma \ref{odomflow}.  Since the net number of particles to enter a site $x\neq o$ is at most one, we have $2d~\div \kappa (x) \geq -1$.  Likewise $2d~\div \kappa(o) = n-1$.  Taking the divergence in (\ref{gradodom}), we obtain
	\begin{align} \Delta u(x) &\leq 1 + \div R(x), \qquad x\neq o;  \label{laplodom}  \\ 
	\Delta u(o) &= 1-n + \div R(o). \label{laplodomat0} 
\end{align}
Let $T$ be the first exit time from $B_r$, and define
	\[ f(x) = \EE_x u(X_T) - \EE_x T + n \EE_x \# \{j<T | X_j=0\}. \]
Then $\Delta f(x) = 1$ for $x \in B_r-\{o\}$ and $\Delta f(o)=1-n$.  Moreover $f \geq 0$ on $\partial B_r$.  It follows from Lemma~\ref{gammaupperbound} with $m=n$ that $f \geq \gamma - C_4$ on $B_r$ for a suitable constant $C_4$.
	
We have
	\[ u(x)-\EE_x u(X_T) = \sum_{k \geq 0} \EE_x \Big( u(X_{k\wedge T}) - u(X_{(k+1)\wedge 
T}) \Big). \] 
Each summand on the right side is zero on the event $\{T \leq k\}$, hence
	\[ \EE_x \Big( u(X_{k\wedge T}) - u(X_{(k+1)\wedge T}) ~|~ \mathcal{F}_{k\wedge 
T} \Big) = -\Delta u (X_k) 1_{\{T>k\}}. \]
Taking expectations and using (\ref{laplodom}) and (\ref{laplodomat0}), we obtain
	\begin{align*}
	u(x)-\EE_x u(X_T) 
	&\geq \sum_{k\geq 0} \EE_x \Big[ 1_{\{T>k\}} (n 1_{\{X_k=o\}} -1 - \div R(X_k) ) \Big] \\
	& \hspace{-0.25in} = n\EE_x \# \{k < T | X_k=o \} - \EE_x T -  \sum_{k\geq 0} \EE_x \Big[ 1_{\{T>k\}}  \div R(X_k) \Big],
	\end{align*}
hence
	\begin{equation} \label{sumofRs}
	u(x) - f(x) \geq - \frac{1}{2d} \sum_{k \geq 0} \EE_x \left[ 1_{\{T>k\}} \sum_{z \sim X_k} R(X_k,z) \right]. \end{equation}
Since random walk exits $B_r$ with probability at least $\frac{1}{2d}$ every time it reaches a site adjacent to the boundary $\partial B_r$, the expected time spent adjacent to the boundary before time $T$ is at most $2d$.  Since $|R|\leq 4d$, the terms in (\ref{sumofRs}) with $z \in \partial B_r$ contribute at most $16d^3$ to the sum.  Thus
	\[ u(x) - f(x) \geq - \frac{1}{2d}  \sum_{k\geq 0} \EE_x  \left[ \sum_{\substack{y,z\in B_r \\ y\sim z}} 1_{\{T>k\} \cap \{X_k=y\}}R(y,z) \right] - 8d^2. \]
For $y \in B_r$ we have $\{X_k = y\} \cap \{T>k\} = \{X_{k\wedge T} = y\}$, hence
	\begin{equation} \label{triplesum} u(x) - f(x) \geq - \frac{1}{2d} \sum_{k\geq 0} \sum_{\substack{y,z \in B_r \\ y\sim z}} \PP_x(X_{k\wedge T}=y)R(y,z) - 8d^2. \end{equation}
Write $p_k(y) = \PP_x(X_{k\wedge T} = y)$.  Note that since $\nabla u$ and $\kappa$ are antisymmetric,~$R$ is antisymmetric.  Thus
	\begin{align*} \sum_{\substack{y,z\in B_r \\ y\sim z}} p_k(y) R(y,z)
	&=  - \sum_{\substack{y,z\in B_r \\ y\sim z}} p_k(z) R(y,z) \\
	&= \sum_{\substack{y,z\in B_r \\ y\sim z}} \frac{p_k(y)-p_k(z)}{2} R(y,z). \end{align*}
Summing over $k$ and using the fact that $|R|\leq 4d$, we conclude from (\ref{triplesum}) that 
	\[ u(x) \geq f(x) - \sum_{\substack{y,z\in B_r \\ y\sim z}} |G(x,y) - G(x,z)| - 8d^2, \]
where $G=G_{B_r}$ is the Green's function for simple random walk stopped on exiting $B_r$.     By Lemma~\ref{greengradbound} we obtain
	\[ u(x) \geq f(x) - C_3 (r+1-|x|) \log \frac{r}{r+1-|x|} - 8d^2. \]
Using the fact that $f \geq \gamma - C_4$, we obtain from Lemma~\ref{gammalowerbound}
	\[ u(x) \geq (r-|x|)^2 - C_3 (r+1-|x|) \log \frac{r}{r+1-|x|} + O\left(\frac{r^d}{|x|^d}\right). \]
The right side is positive provided $r/3 \leq |x| < r- C_5 \log r$.  For $x\in B_{r/3}$, by Lemma~\ref{gammanearorigin} we have $u(x)>r^2/4-C_3 r \log \frac32 >0$, hence $B_{r-C_5\log r} \subset A_n$.
\end{proof}

\subsection{Outer Estimate}

The following result is due to Holroyd and Propp (unpublished); we include a proof for the sake of completeness.  Notice that the bound in (\ref{HPbound}) does not depend on the number of particles.

\begin{prop} \label{HP}
Let $\Gamma$ be a finite connected graph, and let $Y\subset Z$ be subsets of the vertex set of $\Gamma$.  Let $s$ be a nonnegative integer-valued function on the vertices of $\Gamma$.  Let $H_w(s,Y)$ be the expected number of particles stopping in $Y$ if $s(x)$ particles start at each vertex $x$ and perform independent simple random walks stopped on first hitting $Z$.  Let $H_r(s,Y)$ be the number of particles stopping in $Y$ if $s(x)$ particles start at each vertex $x$ and perform rotor-router walks stopped on first hitting $Z$.  Let $H(x) = H_w(1_{x},Y)$.  Then
	\begin{equation} \label{HPbound} 
	| H_r(s,Y) - H_w(s,Y) | \leq  \sum_{u \notin Z} \sum_{v \sim u} | H(u)-H(v) |
	\end{equation}
independent of $s$ and the initial positions of the rotors.
\end{prop} 

\begin{proof}
For each vertex $u \notin Z$, arbitrarily choose a neighbor $\eta(u)$.  Order the neighbors $\eta(u)=v_1, v_2, \ldots, v_d$ of~$u$ so that the rotor at~$u$ points to~$v_{i+1}$ immediately after pointing to $v_i$ (indices mod $d$).  We assign {\it weight} $w(u,\eta(u))=0$ to a rotor pointing from $u$ to $\eta(u)$, and weight $w(u,v_i) = H(u) - H(v_i) + w(u,v_{i-1})$ to a rotor pointing from $u$ to $v_i$.  These assignments are consistent since $H$ is a harmonic function: $\sum_i (H(u) - H(v_i)) = 0$.  We also assign weight $H(u)$ to a particle located at $u$.  The sum of rotor and particle weights in any configuration is invariant under the operation of routing a particle and rotating the corresponding rotor.  Initially, the sum of all particle weights is $H_w(s,Y)$.  After all particles have stopped, the sum of the particle weights is $H_r(s,Y)$.  Their difference is thus at most the change in rotor weights, which is bounded above by the sum in (\ref{HPbound}).
\end{proof}


%
  
For $\rho \in \Z$ let
	\begin{equation} \label{shelldefn} \mathcal{S}_\rho = \{x\in \Z^d ~:~ \rho \leq |x| < \rho+1 \}. \end{equation}
Then
	\[ B_\rho = \{x\in \Z^d ~:~ |x| <\rho \} = \mathcal{S}_0 \cup \ldots \cup \mathcal{S}_{\rho-1}. \]
Note that for simple random walk started in $B_\rho$, the first exit time of $B_\rho$ and first hitting time of $\mathcal{S}_\rho$ coincide.  Our next result is a modification of Lemma~5(b) of \cite{LBG}.

\begin{figure}
\centering
\label{inthezonefigure}
\resizebox{3in}{!}{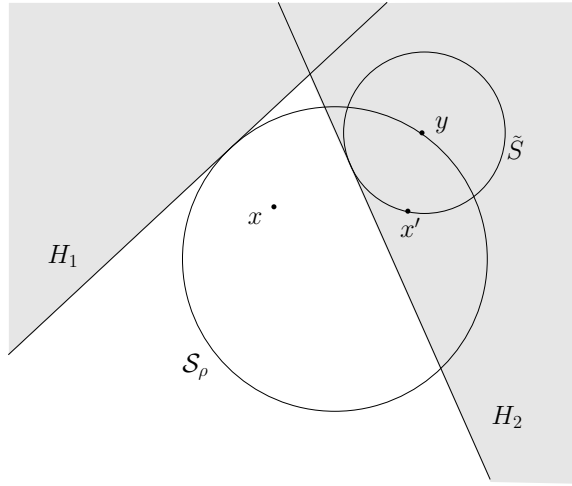}
\caption{Diagram for the proof of Lemma~\ref{inthezone}.}
\end{figure}

\begin{lemma} \label{inthezone}
Fix $\rho \geq 1$ and $y\in \mathcal{S}_{\rho}$.  For $x \in B_\rho$ let $H(x) = \PP_x(X_T=y)$, where $T$ is the first hitting time of $\mathcal{S}_\rho$.  Then
	\begin{equation} \label{spreadingout} H(x) \leq \frac{J}{|x-y|^{d-1}} \end{equation}
for a constant J depending only on $d$.
\end{lemma}

\begin{proof}
We induct on the distance $|x-y|$, assuming the result holds for all $x'$ with $|x'-y| \leq \frac12 |x-y|$; the base case can be made trivial by choosing $J$ sufficiently large.  
By Lemma~5(b) of \cite{LBG}, we can choose $J$ large enough so that the result holds provided $|y|-|x| \geq 2^{-d-3} |x-y|$.  Otherwise, let $H_1$ be the outer half-space tangent to $\mathcal{S}_{\rho}$ at the point of $\mathcal{S}_{\rho}$ closest to $x$, and let $H_2$ be the inner half-space tangent to the ball $\widetilde{S}$ of radius $\frac12 |x-y|$ about $y$, at the point of $\widetilde{S}$ closest to $x$.  By Lemma~\ref{hitthehyperplane} applied to these half-spaces, the probability that random walk started at $x$ reaches $\widetilde{S}$ before hitting $\mathcal{S}_\rho$ is at most $2^{1-d}$.   Writing $\widetilde{T}$ for the first hitting time of $\widetilde{S} \cup \mathcal{S}_{\rho}$, we have
	\[ H(x) \leq \sum_{x' \in \widetilde{S}} \PP_x(X_{\widetilde{T}}=x') H(x') \leq 2^{1-d} J \cdot \left(\frac{|x-y|}{2}\right)^{1-d} \]
where we have used the inductive hypothesis to bound $H(x')$.
\end{proof}
 
The \emph{lazy} random walk in $\Z^d$ stays in place with probability $\frac12$, and moves to each of the $2d$ neighbors with probability $\frac{1}{4d}$.
We will need the following standard result, which can be derived e.g.\ from the estimates in \cite{Lindvall}, section II.12; we include a proof for the sake of completeness. 
 
\begin{lemma}
\label{lazycoupling}
Given $u \sim v \in \Z^d$, lazy random walks started at $u$ and $v$ can be coupled with probability $1-C/R$ before either reaches distance $R$ from $u$, where $C$ depends only on $d$.
\end{lemma}

\begin{proof}
Let $i$ be the coordinate such that $u_i \neq v_i$.  To define a step of the coupling, choose one of the $d$ coordinates uniformly at random.  If the chosen coordinate is different from $i$, let the two walks take the same lazy step so that they still agree in this coordinate.  If the chosen coordinate is $i$, let one walk take a step while the other stays in place.  With probability $\frac12$ the walks will then be coupled.  Otherwise, they are located at points $u',v'$ with $|u'-v'|=2$.  Moreover, $\PP\big(|u-u'|\geq \frac{R}{2\sqrt{d}}\big)<\frac{C'}{R}$ for a constant $C'$ depending only on $d$.  From now on, whenever coordinate $i$ is chosen, let the two walks take lazy steps in opposite directions. 

Let
	\[ H_1 = \left\{x\,\Big|\,x_i = \frac{u'_i+v'_i}{2}\right\} \]
be the hyperplane bisecting the segment $[u',v']$.  Since the steps of one walk are reflections in $H_1$ of the steps of the other, the walks couple when they hit $H_1$.  Let $Q$ be the cube of side length $R/\sqrt{d}+2$ centered at $u$, and let $H_2$ be a hyperplane defining one of the faces of $Q$.  By Lemma~\ref{hitthehyperplane} with $h_1=1$ and $h_2=R/4\sqrt{d}$, the probability that one of the walks exits $Q$ before the walks couple is at most $2d \cdot \frac52 \frac{h_1+1}{h_2} \left(1+\frac{1}{2h_2}\right)^2 \leq 40\, d^{3/2} \big(1+2\sqrt{d}\big)^2 /R$.
\end{proof}

\begin{lemma}
\label{harmonicmeasuregrad}
With $H$ defined as in Lemma~{\em \ref{inthezone}}, we have
	\[ \sum_{u \in B_\rho} \sum_{v \sim u} |H(u)-H(v)| \leq J' \log \rho \]
for a constant $J'$ depending only on $d$.
\end{lemma}

\begin{proof}
Given $u \in B_\rho$ and $v \sim u$, by Lemma~\ref{lazycoupling}, lazy random walks started at $u$ and $v$ can be coupled with probability $1-2C/|u-y|$ before either reaches distance $|u-y|/2$ from $u$.  If the walks reach this distance without coupling, by Lemma~\ref{inthezone} each has still has probability at most $J/|u-y|^{d-1}$ of exiting $B_\rho$ at $y$.  By the strong Markov property it follows that
	\[ |H(u) - H(v)| \leq \frac{2CJ}{|u-y|^d}. \]
Summing in spherical shells about $y$, we obtain
	\[ \sum_{u \in B_\rho} \sum_{v \sim u} |H(u)-H(v)| \leq \sum_{t=1}^{\rho} d\omega_d t^{d-1} \frac{2CJ}{t^d} \leq J' \log \rho. \qed \]
\renewcommand{\qedsymbol}{}
\end{proof}

We remark that Lemma~\ref{harmonicmeasuregrad} could also be inferred from Lemma~\ref{inthezone} using \cite[Thm.\ 1.7.1]{Lawler} in a ball of radius $|u-y|/2$ about $u$.

To prove the outer estimate of Theorem~\ref{rotorcirc}, we will make use of the \emph{abelian property} of rotor-router aggregation.  Fix a finite set $\Gamma \subset \Z^d$ containing the origin.  Starting with $n$ particles at the origin, at each time step, choose a site $x \in \Gamma$ with more than one particle, rotate the rotor at $x$, and move one particle from $x$ to the neighbor the rotor points to.  After a finite number of such choices, each site in $\Gamma$ will have at most one particle, and all particles that exited $\Gamma$ will be on the boundary $\partial \Gamma$.  The abelian property says that \emph{the final configuration of particles and the final configuration of rotors do not depend on the choices}.  For a proof, see \cite[Prop.~4.1]{DF}.

In our application, we will fix $\rho \geq r$ and stop each particle in rotor-router aggregation either when it reaches an unoccupied site or when it reaches the spherical shell $\mathcal{S}_\rho$.  Let $N_\rho$ be the number of particles that reach $\mathcal{S}_\rho$ during this process.  Note that at some sites in $\mathcal{S}_\rho$, more than one particle may have stopped.  If we let each of these extra particles in turn continue performing rotor-router walk, stopping either when it reaches an unoccupied site or when it hits the larger shell $\mathcal{S}_{\rho+h}$, then by the abelian property, the number of particles that reach $\mathcal{S}_{\rho+h}$ will be $N_{\rho+h}$.  We will show that when $h$ is order $r^{1-1/d}$, a constant fraction of the particles that reach $\mathcal{S}_\rho$ find unoccupied sites before reaching $\mathcal{S}_{\rho+h}$.

\begin{proof}[Proof of Theorem~\ref{rotorcirc}, Outer Estimate]
Fix integers $\rho \geq r$ and $h\geq 1$.
In the setting of Proposition~\ref{HP}, let $\Gamma$ be the lattice ball $B_{\rho+h+1}$, and let $Z = \mathcal{S}_{\rho+h}$.  Fix $y \in \mathcal{S}_{\rho+h}$ and let $Y=\{y\}$.  For $x \in \mathcal{S}_\rho$, let $s(x)$ be the number of particles stopped at $x$ if each particle in rotor-router aggregation is stopped either when it reaches an unoccupied site or when it reaches $\mathcal{S}_\rho$.  Write 
	 \[ H(x) = \PP_x(X_T=y) \]
where $T$ is the first hitting time of $\mathcal{S}_{\rho+h}$.   By Lemma~\ref{inthezone} we have
	\begin{equation} \label{H_wbound} H_w(s,y) = \sum_{x \in \mathcal{S}_\rho} s(x)H(x) \leq \frac{J N_\rho}{h^{d-1}} \end{equation}
where
	 \[ N_\rho = \sum_{x \in \mathcal{S}_\rho} s(x) \]
is the number of particles that ever visit the shell $\mathcal{S}_\rho$. 

By Lemma~\ref{harmonicmeasuregrad} the sum in (\ref{HPbound}) is at most $J' \log h$, hence from Propositon~\ref{HP} and (\ref{H_wbound}) we have
	\begin{equation} \label{tinyRRs} H_r(s,y) \leq \frac{J N_\rho}{h^{d-1}} + J' \log h. \end{equation}
Let $\rho(0) = r$, and define $\rho(i)$ inductively by
	\begin{equation} \label{doublemin} \rho(i+1) = \text{min}\left\{\rho(i)+N_{\rho(i)}^{2/(2d-1)},~\text{min} \{\rho> \rho(i) | N_\rho \leq N_{\rho(i)}/2 \}\right\}. \end{equation}
Fixing $h < \rho(i+1)-\rho(i)$, we have
	\[ h^{d-1} \log h \leq N_{\rho(i)}^{\frac{2d-2}{2d-1}} \log N_{\rho(i)} \leq N_{\rho(i)}; \]
so (\ref{tinyRRs}) with $\rho=\rho(i)$ simplifies to
	\begin{equation} \label{eventinierRRs} H_r(s,y) \leq \frac{CN_{\rho(i)}}{h^{d-1}} \end{equation}
where $C = J + J'$.

Since all particles that visit $\mathcal{S}_{\rho(i)+h}$ during rotor-router aggregation must pass through $\mathcal{S}_{\rho(i)}$, we have by the abelian property
	\begin{equation} \label{alltermsaresmall} N_{\rho(i)+h} \leq \sum_{y\in \mathcal{S}_{\rho(i)+h}} H_r(s,y). \end{equation}
Let $M_k = \# (A_n \cap \mathcal{S}_k)$.  There are at most $M_{\rho(i)+h}$ nonzero terms in the sum on the right side of (\ref{alltermsaresmall}), and each term is bounded above by (\ref{eventinierRRs}), hence
	\[ M_{\rho(i)+h} \geq N_{\rho(i)+h} \frac{h^{d-1}}{CN_{\rho(i)}} \geq  \frac{h^{d-1}}{2C} \]
where the second inequality follows from $N_{\rho(i)+h} \geq N_{\rho(i)}/2$.  Summing over $h$, we obtain
	\begin{equation} \label{lotsavolume}  \sum_{\rho=\rho(i)+1}^{\rho(i+1)-1} M_\rho \geq \frac{1}{2dC} (\rho(i+1)-\rho(i)-1)^d. \end{equation}
The left side is at most $N_{\rho(i)}$, hence
	\[ \rho(i+1)-\rho(i) \leq (2dC N_{\rho(i)})^{1/d} \leq N_{\rho(i)}^{2/(2d-1)} \]
provided $N_{\rho(i)} \geq C' := (2dC)^{2d-1}$.
Thus the minimum in (\ref{doublemin}) is not attained by its first argument.  It follows that $N_{\rho(i+1)}\leq N_{\rho(i)}/2$, hence $N_{\rho(a \log r)} < C' $ for a sufficiently large constant $a$.

By the inner estimate, since the ball $B_{r-c\log r}$ is entirely occupied, we have
	\begin{align*} \sum_{\rho \geq r} M_\rho &\leq \omega_d r^d - \omega_d (r-c\log r)^d \\
	&\leq cd\omega_d r^{d-1} \log r. \end{align*}
Write $x_i = \rho(i+1)-\rho(i)-1$; by (\ref{lotsavolume}) we have
	\[ \sum_{i=0}^{a \log r} x_i^d \leq cd\omega_d r^{d-1} \log r, \]
By Jensen's inequality, subject to  this constraint, $\sum x_i$ is maximized when all $x_i$ are equal, in which case $x_i \leq C'' r^{1-1/d}$ and 
	\begin{equation} \rho(a \log r) = r+\sum x_i \leq r+C''r^{1-1/d} \log r. \label{outershellbound} \end{equation}
Since $N_{\rho(a \log r)} < C'$ we have $N_{\rho(a \log r)+C'} = 0$; that is, no particles reach the shell $\mathcal{S}_{\rho(a \log r)+C'}$.  Taking $c' = C' + C''$, we obtain from (\ref{outershellbound})
	\[ A_n \subset B_{r (1 + c'r^{-1/d}\log r)}. \qed \]
\renewcommand{\qedsymbol}{}
\end{proof}

\section{Concluding Remarks}

\begin{figure}
\label{invertedrotor}
\centering
\includegraphics{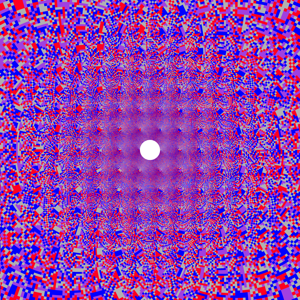}
\caption{Image of the rotor-router aggregate of one million particles under the map $z \mapsto 1/z^2$.  The colors represent the rotor directions.  The white disc in the center is the image of the complement of the occupied region.}
\end{figure}

A number of intriguing questions remain unanswered.  Although we have shown that the asymptotic shape of the rotor-router model is a ball, the near perfect circularity found in Figure~\ref{rotor1m} remains a mystery.  In particular, we do not know whether an analogue of Theorem~\ref{divsandcircintro} holds for the rotor-router model, with constant error in the radius as the number of particles grows.   

Equally mysterious are the patterns in the rotor directions evident in Figure~\ref{rotor1m}.  The rotor directions can be viewed as values of the odometer function mod $2d$, but our control of the odometer is not fine enough to provide useful information about the patterns.  If the rescaled occupied region $\sqrt{\pi/n} \,A_n$ is viewed as a subset of the complex plane, it appears that the monochromatic regions visible in Figure~\ref{rotor1m}, in which all rotors point in the same direction, occur near points of the form $(1+2z)^{-1/2}$, where $z = a+bi$ is a Gaussian integer (i.e.\ $a,b \in \Z$).  We do not even have a heuristic explanation for this phenomenon.  Figure~7 shows the image of $A_{1,000,000}$ under the map $z \mapsto 1/z^2$; the monochromatic patches in the transformed region occur at lattice points.

L\'{a}szl\'{o} Lov\'{a}sz (personal communication) has asked whether the occupied region $A_n$ is simply connected, i.e.\ whether its complement is connected.  While Theorem~\ref{rotorcirc} shows that $A_n$ cannot have any holes far from the boundary, we cannot answer his question at present.

A final question is whether our methods could be adapted to internal DLA to show that if $n=\omega_d r^d$, then with high probability $B_{r-c\log r} \subset I_n$, where $I_n$ is the internal DLA cluster of $n$ particles.  The current best bound is due to Lawler \cite{Lawler95}, who proves that with high probability $B_{r-r^{1/3}(\log r)^2} \subset I_n$.

\section*{Acknowledgments}

The authors thank Jim Propp for bringing the rotor-router model to our attention and for many fruitful discussions.  We also had useful discussions with Oded Schramm, Scott Sheffield and Misha Sodin.  We thank Wilfried Huss for pointing out an error in an earlier draft.  Yelena Shvets helped draw some of the figures.








\end{document}